\documentclass[12pt, reqno]{amsart}
\usepackage{amscd,amsmath,amsthm,amssymb,graphics}
\usepackage{amsfonts,amssymb,amscd,amsmath,enumitem,verbatim}
\usepackage[a4paper,top=3cm,left=3cm,right=3cm]{geometry}
\theoremstyle{plain}
\usepackage{booktabs}
\usepackage{hyperref}
\newtheorem{Theorem}{Theorem}

\newtheorem{Proposition}[Theorem]{Proposition}

\newtheorem{Conjecture}[Theorem]{Conjecture}

\theoremstyle{definition}
\newtheorem{Definition}[Theorem]{Definition}
\newtheorem{Remark}[Theorem]{Remark}

\newtheorem{Example}[Theorem]{Example}

\newcommand\blfootnote[1]{%
  \begingroup
  \renewcommand\thefootnote{}\footnote{#1}%
  \addtocounter{footnote}{-1}%
  \endgroup
}

\title{Continued fractions in the field of p-adic numbers}
\author{Giuliano Romeo}

\begin{document}

\maketitle

\begin{center}
Department of Mathematical Sciences Giuseppe Luigi Lagrange,\\
Politecnico of Torino, Corso Duca degli Abruzzi 24, Torino
\ \\ \ \\
giuliano.romeo@polito.it, 
\end{center}

\begin{abstract}
Continued fractions have a long history in number theory, especially in the area of Diophantine approximation. The aim of this expository paper is to survey the main results on the theory of $p$--adic continued fractions, i.e. continued fractions defined over the field of $p$--adic numbers $\mathbb{Q}_p$, which in the last years has recorded a considerable increase of interest and research activity. We start from the very first definitions up to the most recent developments and open problems.
\end{abstract}

\tableofcontents

\section{Introduction}
Continued fractions are objects of the form \begin{equation}
a_0+\cfrac{b_1}{a_1+\cfrac{b_2}{a_2+\ddots}},
\end{equation}
and they have a long history in number theory, especially in the area of Diophantine approximation. They have been extensively explored over the real numbers, where they provide finite representations for rational numbers and periodic representations for quadratic irrationals by means of integer sequences (Lagrange's Theorem). Moreover, the classical algorithm for continued fractions in $\mathbb{R}$ provides, at each step, the best rational approximation of irrational numbers. Because of their optimal properties, they have been employed in various areas of mathematics and there exist several generalizations of the classical theory of continued fractions. In 1940, Mahler \cite{MAH} raised the problem of defining a continued fraction algorithm working in the field of $p$--adic numbers $\mathbb{Q}_p$. The field of $p$--adic numbers is largely studied in algebraic number theory and it is obtained as the completion of the field of rational numbers $\mathbb{Q}$ with respect to the $p$--adic absolute value. The completion is the smallest field containing $\mathbb{Q}$ in which all Cauchy sequences with respect to the $p$--adic absolute value are convergent. The problem of defining a $p$--adic continued fraction algorithm sharing all the optimal properties enjoyed by classical continued fractions is still open. In particular, an algorithm providing periodic representations for every $p$--adic quadratic irrational is not known, i.e. a $p$--adic analogue of Lagrange's Theorem does not exist yet. The main problem is that there is not an intuitive satisfying definition for the integer part of a $p$--adic number. There are at least two very natural definitions, due to Schneider \cite{SCH} and Ruban \cite{RUB}, both around 1970. The two algorithms try to replicate the standard algorithm for real continued fractions in two different senses, that we examine in Section \ref{Sec: algo}. Historically, another of the main algorithms for $p$--adic continued fractions is due to Browkin \cite{BI}. His algorithm is very similar to Ruban's one, but it has the very important property of producing finite continued fractions for all rational numbers. In fact, Ruban's and Schneider's continued fractions are not always finite over $\mathbb{Q}$. Moreover, Schneider's and Ruban's algorithms have been proved to be not periodic for every quadratic irrational \cite{TIL,CVZ}. Although a proof has not been provided, Lagrange's Theorem seems to fail also for Browkin's algorithm (see the numerical simulations in \cite{MR}). For these reasons, the research of new definitions for expressing a $p$--adic number as a continued fraction has been widely developed, particularly in the last few years. Together with the design of new algorithms, the analysis of the properties of the existent algorithms has been developed in several directions, inspired by the numerous fields in which classical continued fractions have been employed throughout the centuries. \\

The purpose of this survey is to give a concise overview of the main results, developments and open problems in the theory of $p$--adic continued fractions. In Section \ref{Sec: preli}, we recall some notation for the classical theory of continued fractions and $p$--adic numbers. In Section \ref{Sec: algo} we describe the various algorithms that have been defined over the years, underlining their motivations and their main properties. In Section \ref{Sec: convergence}, we discuss the $p$--adic convergence of a continued fraction, which is the very first requirement for the definition of an algorithm. In Section \ref{Sec: ratio}, we present the properties of the expansions of rational numbers that, as already mentioned, are not always finite. Section \ref{Sec: periodicity} contains all the main results related to the periodicity of $p$--adic continued fractions and the main developments towards a $p$--adic analogue of Lagrange's Theorem. In Section \ref{Sec: approxi}, we analyze the properties of approximation of continued fractions with respect to the $p$--adic absolute value. In Section \ref{Sec: transce}, we collect some recent results on the transcendence of $p$--adic continued fractions, inspired by many famous transcendence criteria for classical continued fractions. Finally, in Section \ref{Sec: multidim} we introduce the generalization of $p$--adic continued fractions to higher dimensions, in analogy with Jacobi-Perron algorithm for real multidimensional continued fractions.

\section{Some notation}\label{Sec: preli}
In this section we recall some useful results from the theory of continued fractions and the field of $p$--adic numbers and we fix some notation that we are going to use throughout the survey. For more background on continued fractions and $p$--adic numbers we refer the reader, respectively, to \cite{OLDS, WALL} and \cite{FG}.\\

We call \textit{continued fraction} an object of the form
\begin{equation}
a_0+\cfrac{b_1}{a_1+\cfrac{b_2}{a_2+\ddots}},
\end{equation}
and we denote it by
\[\begin{bmatrix}
\ b_1 & b_2 & \ldots  \  \\
\ a_0 & a_1 & a_2 & \ldots  \ \\
\end{bmatrix},\]
where the coefficients $a_i$ and $b_i$ are elements in a field and the expansion can be either finite or infinite. If $b_i=1$ for all $i$, the continued fraction is called \textit{simple}, that is
\[a_0+\cfrac{1}{a_1+\cfrac{1}{a_2+\ddots}},\]
and we denote it by $[a_0,a_1,a_2,\ldots]$. The coefficients $a_0,a_1,a_2\ldots$ are called \textit{partial quotients}. For all $n\in\mathbb{N}$, the rational number
\[\frac{A_n}{B_n}=\begin{bmatrix}
\ b_1 & \ldots & b_{n-1} & b_n & \  \\
\ a_0 & a_1 & \ldots & a_{n-1} & a_n  \ \\
\end{bmatrix}=
a_0+\cfrac{b_1}{a_1+\ddots + \cfrac{b_{n-1}}{a_{n-1}+\cfrac{b_n}{a_n}}},\]
corresponding to the continued fraction stopped at the $n$-th term, is called $n$-th \textit{convergent} of the continued fraction. The $A_n$'s are called \textit{partial numerators} and the $B_n$'s are called \textit{partial denominators}. The sequences $\{A_n\}_{n\in\mathbb{N}}$ and $\{B_n\}_{n\in\mathbb{N}}$ satisfy the recursions
\begin{equation}\label{Eq: Recursions}
\begin{cases}
A_0=a_0,\\
A_1=a_1a_0+b_1,\\
A_n=a_nA_{n-1}+b_nA_{n-2}, \ \ n \geq 2,
\end{cases}
\begin{cases}
B_0=1,\\
B_1=a_1,\\
B_n=a_nB_{n-1}+b_nB_{n-2}, \ \ n \geq 2.
\end{cases}    
\end{equation}
The partial numerators and denominators of the convergents can be represented also using the following matrix form. For all $n\in\mathbb{N}$,
\begin{equation}\label{Eq: matrix}
\begin{pmatrix}
\ A_n & A_{n-1} \ \\
\ B_n & B_{n-1}  \ 
\end{pmatrix}
=
\begin{pmatrix}
\ a_0 & 1 \ \\
\ 1 & 0  \ 
\end{pmatrix}
\begin{pmatrix}
\ 1 & 0 \ \\
\ 0 & b_1  \ 
\end{pmatrix}
\begin{pmatrix}
\ a_1 & 1 \ \\
\ 1 & 0  \ 
\end{pmatrix}
\ldots
\begin{pmatrix}
\ 1 & 0 \ \\
\ 0 & b_n  \ 
\end{pmatrix}
\begin{pmatrix}
\ a_n & 1 \ \\
\ 1 & 0  \ 
\end{pmatrix}.
\end{equation}

Furthermore, for an infinite continued fraction representing an element $\alpha$, we introduce the sequence of \textit{complete quotients} $\{\alpha_n\}_{n\in\mathbb{N}}$ defined, for all $n\in\mathbb{N}$, as
\[\alpha_{n+1}=\frac{b_n}{\alpha_n-a_n},\]
starting from $\alpha_0=\alpha$. We say that an infinite continued fraction is \textit{periodic} (or \textit{eventually periodic}) if and only if there exist $n_0\in\mathbb{N}$ and $k\geq 1$ such that $a_{n+k}=a_n$ $b_{n+k}=b_n$ for all $n\geq n_0$. In this case we call $k$ the \textit{period length} and $n_0$ the \textit{pre-period length} (or just \textit{period} and \textit{pre-period} where there is no risk of ambiguity).
If $n_0=0$, the continued fraction is said \textit{purely periodic} and the period starts from the begin, without a pre-period. The standard algorithm to express a real number $\alpha$ through a simple continued fraction $[a_0,a_1,\ldots]$ is the following, with $\alpha_0=\alpha$:
\begin{equation}\label{Alg: R}
\begin{cases}
a_n=\lfloor \alpha_n \rfloor \\
\alpha_{n+1}=\frac{1}{\alpha_n-a_n},
\end{cases}
\end{equation}

where $\lfloor \alpha_n \rfloor $ denotes the integer part of $\alpha$, that is the greatest integer $a_n\leq \alpha_n$. If $\alpha_n=a_n$ for some $n\in\mathbb{N}$, then the algorithm terminates and the continued fraction is finite.\\

In the following, let $p$ be an odd prime number and let us denote by $v_p(\cdot)$ the \textit{$p$--adic valuation} and by $|\cdot|$ and $|\cdot|_p$, respectively, the \textit{Euclidean absolute value} and the \textit{$p$--adic absolute value}. The field of $p$--adic numbers $\mathbb{Q}_p$ is the completion of $\mathbb{Q}$ with respect to the $p$--adic absolute value, i.e. the smallest field containing $\mathbb{Q}$ in which all the Cauchy sequences are convergent. The field of $p$--adic numbers contains all the power series in $p$ with finite tail, that is:

\begin{equation}\label{Eq: padicseries}
\mathbb{Q}_p=\left\{ \sum\limits_{n=-r}^{+\infty} c_np^n \ \Big| \ r\in\mathbb{Z}, \ c_n\in\mathbb{Z}/p\mathbb{Z} \right\}.
\end{equation}

By construction $\mathbb{Q}\subset\mathbb{Q}_p$ and it is a known fact that rational numbers correspond exactly to the set finite and eventually periodic series. Moreover, we call $\mathbb{Z}_p$ the set of \textit{$p$--adic integers},
\[\mathbb{Z}_{p}=\Bigg\{ \sum\limits_{n=0}^{+\infty} c_np^n \ \Big| \ c_n\in \mathbb{Z}/p\mathbb{Z} \Bigg\},\]
that is the set of all $p$--adic numbers with non-negative valuation.

\section{Algorithms for $p$--adic continued fractions}\label{Sec: algo}
In this section we present the main definitions for continued fractions in the field of $p$--adic numbers that have appeared in literature. The first attempt is to emulate Algorithm  the standard algorithm \eqref{Alg: R} for real continued fractions. However, the definition of integer part of a $p$--adic number $\alpha\in\mathbb{Q}_p$ is not uniquely determined. In $\mathbb{R}$, the integer part of a real number $\alpha$ is defined as the unique integer $b$ less than $\alpha$ for which $|\alpha-a|<1$. However, in $\mathbb{Q}_p$ there are infinitely many integers such that this happens, since there are infinitely many $a\in\mathbb{Z}$ such that $p$ divides $\alpha-a$. One of the first definitions of $p$-adic continued fractions is due to Schneider \cite{SCH}. Let us observe that all the integers $a$ satisfying $0\leq |\alpha -a|_p<1$ are congruent modulo $p$. Therefore, it is meaningful to define the partial quotients as the unique representatives $a_n$ of $a$ that lie inside $\{0,\ldots,p-1\}$. Schneider followed this approach to provide a non-simple continued fraction expansion
 \begin{equation}
a_0+\cfrac{b_1}{a_1+\cfrac{b_2}{a_2+\ddots}},
\end{equation}
for all $p$-adic integers. For $\alpha_0\in \mathbb{Z}_{p}$, with representation
\[\alpha_0=\sum\limits_{n=0}^{+\infty}c_np^n, \ \ \ \ \ c_n\in\{0,\ldots,p-1\},\]
Schneider's algorithm works as follows. The first partial quotient is $a_0=c_0$. Then, for all $n\in\mathbb{N}$,
\begin{align} 
\begin{cases}\label{Alg: Sch}
e_{n+1}=v_p(\alpha_{n}-a_n)\\
b_{n+1}=p^{e_{n+1}}\\
\alpha_{n+1}=\dfrac{b_{n+1}}{\alpha_n-a_n}\\
a_{n+1}=c_0^{(n+1)},\\
\end{cases}
\end{align}
where $c_0^{(n+1)}$ denotes the first term of the expansion of
\[p^e\alpha_{n+1}=\sum\limits_{i=0}^{+\infty}c_i^{(n+1)}p^i.\]
\begin{Example}
Let us give an idea on how Schneider's algorithm \eqref{Alg: Sch} works to provide the $5$--adic continued fraction of $\frac{2}{7}$. In $\mathbb{Q}_5$,
\[\frac{2}{7}=1+2\cdot5+5^2+\ldots,\]
hence $a_0=2$ and $e_1=v_p(\alpha_0-a_0)=v_p(2\cdot5+5^2+\ldots)=1$.
The next complete quotient is then
\[\alpha_1=\frac{p^{e_1}}{\alpha_0-a_0}=\frac{5}{2\cdot5+5^2+\ldots}=\frac{1}{2+5+\ldots}=3+3\cdot5+4\cdot 5^2+\ldots.\]
Hence, $a_1=3$ and the algorithm goes on as before. The expansion becomes eventually periodic with $b_n=5$ for all $n$ and $[a_0,a_1,\ldots]=[3,2,3,\overline{4}]$.
\end{Example}

Ruban's \cite{RUB} and Browkin's \cite{BI, BII} approach is to construct a simple continued fraction by defining an intuitive $p$--adic analogue of the integer part in $\mathbb{R}$.  For any $p$-adic number $\alpha=\sum\limits_{n=-r}^{+\infty} c_np^n\in\mathbb{Q}_p$, $c_n\in\{0,\ldots,p-1\}$,
the integer part used by Ruban is
\[\lfloor \alpha \rfloor_p =\sum\limits_{n=-r}^{0} c_np^n,\]
and $\lfloor \alpha \rfloor_p=0$ if $r<0$. Notice that in this case $\lfloor \alpha \rfloor_p$ is, in general, a rational number. Ruban's continued fractions are simple and the coefficients of the expansion can be computed iteratively by the following algorithm, starting with $\alpha_0=\alpha$.
\begin{align} 
\begin{cases}\label{Alg: Rub}
a_n=\lfloor \alpha_n \rfloor_p\\
\alpha_{n+1}=\frac{1}{\alpha_n-a_n}.
\end{cases}
\end{align}
If at some point $\alpha_n=a_n$, then the algorithm stops and $\alpha=[a_0,\ldots,a_n]$, i.e. $\alpha$ has finite Ruban's continued fraction. The choice of this integer part is, somehow, natural. In fact, in the series $\sum\limits_{n=-r}^{+\infty} c_np^n$, the positive powers of $p$ have fractional $p$--adic absolute value and the floor function $\lfloor \cdot \rfloor_p$ takes the part that has integral absolute value.
\begin{Example}\label{Exa: Ruban}
Let us compute Ruban's expansion of $\alpha_0=-\frac{2}{5}$ in $\mathbb{Q}_7$. The $7$--adic expansion of $\alpha_0$ is
\[\alpha_0 = 1 + 4 \cdot 7 + 5\cdot 7^2 + \ldots,\]
so that $a_0= \lfloor\alpha_0 \rfloor_p = 1$ and
\[\alpha_1 =\frac{1}{\alpha_0-a_0}=\frac{1}{-\frac{2}{5}-1}= -\frac{5}{7} = 2 \cdot 7^{-1} +6 +6 \cdot 7 + \ldots.\]
The second partial quotient is then
\[a_1=\lfloor\alpha_1 \rfloor_p=\frac{2}{7}+6=\frac{44}{7},\]
and we compute the third complete quotient as
\[\alpha_2 =\frac{1}{\alpha_1-a_1}=\frac{1}{-\frac{5}{7}-\frac{44}{7}}= -\frac{1}{7} = 6 \cdot 7^{-1} +6 +6 \cdot 7 + \ldots.\]
Then
\[a_2=\lfloor \alpha_2 \rfloor_p=\frac{48}{7},\]
and we find
\[\alpha_3=\frac{1}{\alpha_2-a_2}=\frac{1}{-\frac{1}{7}-\frac{48}{7}}=-\frac{1}{7}=\alpha_2.\]
 We have found a repetition on the complete quotients, $\alpha_3=\alpha_2$, so the continued fraction repeats from this point onward. It means that
\[-\frac{2}{5}=\left[1,\frac{44}{7},\overline{\frac{48}{7}} \right],\]
hence $-\frac{2}{5}$ has a periodic Ruban's continued fraction in $\mathbb{Q}_7$.
\end{Example}

Browkin's first algorithm  \cite{BI} is very similar to Ruban's algorithm, with the exception that the representatives of $\mathbb{Z}/p\mathbb{Z}$ are chosen in $\{-\frac{p-1}{2},\ldots,\frac{p-1}{2}\}$ instead of $\{0,\ldots,p-1\}$. This small variation is fundamental since Browkin's algorithm produces a finite continued fraction for each rational number, while Ruban's algorithm can also be periodic over the rationals (more details are given in Section \ref{Sec: ratio}). Given $\alpha=\sum\limits_{i=-r}^{+\infty}c_ip^i\in\mathbb{Q}_p$, with $c_i \in \{ -\frac{p-1}{2}, \ldots, \frac{p-1}{2} \}$, Browkin defines the floor function $s:\mathbb{Q}_p\rightarrow \mathbb{Q}$ as
\begin{equation}\label{Eq: sfunc}
s(\alpha)=\sum\limits_{i=-r}^{0}c_ip^i,
\end{equation}
and $s(\alpha)=0$ if $r<0$. Browkin's first algorithm, that we call \textit{Browkin I}, work as follows. At the first step $\alpha_{0}=\alpha$ and, for all $n\geq 0$,
\begin{align} 
\begin{cases}\label{Alg: Br1}
a_n=s(\alpha_n)\\
\alpha_{n+1}=\frac{1}{\alpha_n-a_n}.
\end{cases}
\end{align}
If at some point $\alpha_n=a_n$, then the algorithm stops and $\alpha=[a_0,\ldots,a_n]$, i.e. $\alpha$ has finite \textit{Browkin I} continued fraction. 

\begin{Example}
Let us compute \textit{Browkin I} continued fraction of $\alpha_0=-\frac{2}{5}$ in $\mathbb{Q}_7$ and let us compare with Ruban's expansion obtained in Example \ref{Exa: Ruban}. For Browkin's algorithm, we have to consider the $7$--adic expansion with representatives in $\{-3,\ldots,3\}$, hence:
\[\alpha_0 = 1 -3 \cdot 7 -  7^2 + \ldots.\]
Also in this case the first partial quotient is $a_0= s(\alpha_0 )= 1$ and then
\[\alpha_1 =\frac{1}{\alpha_0-a_0}=\frac{1}{-\frac{2}{5}-1}= -\frac{5}{7} = 2 \cdot 7^{-1} -1 \ldots.\]
The second partial quotient is then
\[a_1=s(\alpha_1)=\frac{2}{7}-1=-\frac{5}{7}=\alpha_1.\]
Therefore,
\[-\frac{2}{5}=\left[1,-\frac{5}{7}\right],\]
hence $-\frac{2}{5}$ has a finite \textit{Browkin I} continued fraction. The difference with Ruban's expansion of Example \ref{Exa: Ruban}, which is periodic, relies only on the choice of the representatives in $\{-3,\ldots,3\}$ instead of $\{0,\ldots,6\}$.
\end{Example}

In 2000, more than 20 years after the first algorithm, Browkin \cite{BII} defined another floor function, that is similar to the first function $s$, but excluding the constant term. For $\alpha=\sum\limits_{i=-r}^{+\infty}c_ip^i\in\mathbb{Q}_p$, with $c_i \in \{ -\frac{p-1}{2}, \ldots, \frac{p-1}{2} \}$, the second floor function is the function $t:\mathbb{Q}_p\rightarrow \mathbb{Q}$, such that
\[t(\alpha)=\sum\limits_{i=-r}^{-1}c_ip^i,\]
and $t(\alpha)=0$ if $r\leq 0$. The second algorithm, that we call \textit{Browkin II}, works on an input $\alpha$ as follows. At the first step $\alpha_{0}=\alpha$ and, for all $n\geq 0$,
\begin{align} 
\begin{cases}\label{Alg: Br2}
a_n=s(\alpha_n) \ \ \ \ \ & \textup{if} \ n \ \textup{even}\\
a_n=t(\alpha_n) & \textup{if} \ n \ \textup{odd}\ \textup{and} \ v_p(\alpha_n-t(\alpha_n))= 0\\
a_n=t(\alpha_n)-sign(t(\alpha_n)) & \textup{if} \ n \ \textup{odd} \ \textup{and} \ v_p(\alpha_n-t(\alpha_n))\neq 0\\
\alpha_{n+1}=\frac{1}{\alpha_n-a_n}.
\end{cases}
\end{align}
If at some point $\alpha_n=a_n$, then the algorithm stops and $\alpha=[a_0,\ldots,a_n]$, i.e. $\alpha$ has finite \textit{Browkin II} continued fraction. The result of the alternation is that, in \textit{Browkin II}, all the even partial quotients are integers and all the odd partial quotients are rationals. The use of the sign function is due to the convergence condition of Proposition \ref{Prop: ConvBr2} proved by Browkin in \cite{BII}.

\begin{Example}
Let us consider $\alpha_0 = \frac{22}{7} \in \mathbb Q_5$ and let us compute its \textit{Browkin II} expansion. The $5$--adic expansion of $\alpha_0$ is
\[\alpha_0 = 1 - 1 \cdot 5 + 1 \cdot 5^2 + \ldots,\]
so that $a_0= s(\alpha_0) = 1$ and
\[\alpha_1 =\frac{1}{\alpha_0-a_0}=\frac{1}{\frac{22}{7}-1}= \frac{7}{15} = -1 \cdot 5^{-1} - 1 +2 \cdot 5 + \ldots.\]
Now we apply the function $t$, thus obtaining $b_1 = t(\alpha_1) = -\frac{1}{5}$ and
\[\alpha_2 =\frac{1}{\alpha_1-a_1}=\frac{3}{2}=-1 -2 \cdot 5 - 2 \cdot 5^2+ \ldots.\]
At the next step we have $a_2 = s(\alpha_2) = -1$ and  $\alpha_3 = \frac{2}{5}$. At this point, since $v_p(\alpha_3-t(\alpha_3))>0$, then we use the sign function, obtaining:
\[a_3 = t(\alpha_3) - sign(t(\alpha_3)) =\frac{2}{5}-1 = -\frac{3}{5}.\]
Therefore $\alpha_4=s(\alpha_4)=b_4=1$ and the expansion of is $\frac{22}{7} = \left[ 1, -\frac{1}{5}, -1, -\frac{3}{5}, 1 \right]$.
\end{Example}
It has been experimentally observed that \textit{Browkin II} produces more periodic expansions for quadratic irrationals than \textit{Browkin I}, hence getting closer to a $p$--adic analogue of Lagrange's theorem (see, for example, \cite{MR}). For this reason, \textit{Browkin II} has been taken as a starting point for the definition of some new algorithms with the aim of improving furthermore its properties of periodicity. In \cite{BCMII}, it is studied a a variant of \textit{Browkin II} in which the representatives are taken in $\{0,\ldots,p-1 \}$ instead of $\{-\frac{p-1}{2},\ldots,\frac{p-1}{2}\}$. This choice improves the properties of approximation \textit{Browkin II} and, unexpectedly, does not compromise the finiteness of the continued fraction expansions of rational numbers. In \cite{MRSII}, it is proposed a $3$-steps generalization of \textit{Browkin II} relying on the convergence condition of Theorem \ref{Thm: ConvBr3} in Section \ref{Sec: convergence}. In order to satisfy that convergence condition, a new integer part is defined, that acts on a $p$--adic integer with zero valuation and leaves apart another $p$--adic integer with zero valuation. For $\alpha=\sum\limits_{n=0}^{+\infty}c_ip^i\in\mathbb{Q}_p$, with $c_i \in \{ -\frac{p-1}{2}, \ldots, \frac{p-1}{2} \}$, the new integer part is
\begin{align}\label{Eq: u}
u(\alpha)=\begin{cases} +1 \ \  &\textup{if}  \ c_0\in\Big\{+2,\ldots,\dfrac{p-1}{2}\Big\}\cup \{-1\}\\
-1 &\textup{if}  \ c_0\in\Big\{-\dfrac{p-1}{2},\ldots,-2\Big\}\cup \{+1\}\\
0 \ \ &\textup{if}  \ c_0=0.
\end{cases}
\end{align}
In \cite{MRSII}, two new algorithms have been defined by exploiting the three integer parts $s,t,u$ in three different steps. For one of them, the authors proved that all rational numbers have a finite continued fraction, hence becoming more interesting than the other one. On input $\alpha_0\in\mathbb{Q}$ the algorithm works as follows, for all $n\in\mathbb{N}$:

\begin{align}\label{Alg: Br3}\begin{cases}
a_n=s(\alpha_n) \ \ \ \ \ &\textup{if} \ n \equiv 0\bmod 3\\
a_n=t(\alpha_n) &\textup{if}  \ n \equiv 1 \bmod 3 \ \textup{and} \ v_p(\alpha_n-t(\alpha_n))= 0\\
a_n=t(\alpha_n)-sign(t(\alpha_n)) & \textup{if} \ n \equiv 1 \bmod 3 \ \textup{and} \  v_p(\alpha_n-t(\alpha_n))\neq0\\
a_n=s(\alpha_n)-u(\alpha_n)  & \textup{if} \ n \equiv 2 \bmod 3\\
\alpha_{n+1}=\frac{1}{\alpha_n-a_n}.
\end{cases}
\end{align}
If at some point $\alpha_n=a_n$, then the algorithm stops and $\alpha=[a_0,\ldots,a_n]$, i.e. $\alpha$ has a finite continued fraction. Finally, in \cite{MR}, the authors proposed a variant of \textit{Browkin II} without the use of the sign function. For all $\alpha_0\in\mathbb{Q}_p$, their algorithm works as follows:
\begin{align} 
\begin{cases}\label{Alg: new}
a_n=s(\alpha_n) \ \ \ \ \ & \textup{if} \ n \ \textup{even}\\
a_n=t(\alpha_n) & \textup{if} \ n \ \textup{odd}\\
\alpha_{n+1}=\frac{1}{\alpha_n-a_n}.
\end{cases}
\end{align}
If at some point $\alpha_n=a_n$, then the algorithm stops and $\alpha=[a_0,\ldots,a_n]$, i.e. $\alpha$ has a finite continued fraction. This choice turns out to improve the periodicity properties of Browkin's second algorithm (more details are given in Section \ref{Sec: periodicity}).\\

However, although all these algorithms improve on several aspects the known methods to generate $p$--adic continued fractions, they experimentally seem still far from an analogue of Lagrange's Theorem and, in general, from a satisfactory algorithm reproducing the same properties of classical continued fractions in $\mathbb{R}$.\\

Very recently, in \cite{CMT2}, the definition of $p$--adic continued fractions has been generalized also to number fields, addressing some questions similar to Rosen \cite{ROS1,ROS2} in the archimedean setting. The authors gave a general definition of $\mathfrak{P}$--adic continued fractions for a prime ideal $\mathfrak{P}$ of the ring of integers $\mathcal{O}_K$ of a number field $K$. Moreover, they investigated their finiteness and periodicity properties, focusing on a general number field $K$ and obtaining some more effective results for the quadratic case. In \cite{MULA}, the construction of \cite{CMT2} has been generalized to quaternion algebras over $\mathbb{Q}$ ramified at $p$, hence including all the classical $p$--adic framework.

\section{Convergence in $\mathbb{Q}_p$}\label{Sec: convergence}
In the field of real numbers, every continued fraction with positive partial quotients converges to an $\alpha\in\mathbb{R}$. This is not always the case in $\mathbb{Q}_p$. The field of $p$--adic numbers is complete with respect to the $p$--adic absolute value $|\cdot|_p$. Therefore, a sequence is convergent in $\mathbb{Q}_p$ if and only if it is a Cauchy sequence. Then, a continued fraction
\[a_0+\cfrac{b_1}{a_1+\cfrac{b_2}{a_2+\ddots}},\]
converges to a $p$--adic number if and only if the sequence of the convergents $\{\frac{A_n}{B_n}\}_{n\in\mathbb{N}}$ is a Cauchy sequence with respect to $|\cdot|_p$. For a non-archimedean absolute value this is equivalent to require that
\[ \lim \limits_{n \rightarrow +\infty} \left|\frac{A_{n+1}}{B_{n+1}}-\frac{A_{n}}{B_{n}} \right|_p=0.\]
The latter quantity can be written as
\begin{equation}\label{Eq: conve}
\left|\frac{A_{n+1}}{B_{n+1}}-\frac{A_{n}}{B_{n}} \right|_p=\left|\frac{(-1)^n}{B_{n}B_{n+1}}\right|_p=p^{v_p(B_nB_{n+1})},
\end{equation}
so that the continued fraction $[a_0,a_1,\ldots]$ converges to an element of $\mathbb{Q}_p$ if and only if
\[\lim\limits_{n\rightarrow +\infty} v_p(B_nB_{n+1})=-\infty.\]
This is the minimum requirement in order to define meaningful continued fractions in the field of $p$--adic numbers, since otherwise a continued fraction does not represent any element of $\mathbb{Q}_p$. A proof of the convergence of Schneider's continued fractions can be found in \cite{VDP} and it exploits the matrix relations \eqref{Eq: matrix}.

\begin{Remark}\label{Rem: ConveBr}
In \cite{BI}, Browkin proved that, for \textit{Browkin I},
\begin{equation}
\begin{split}
v_p(A_n)&=v_p(a_0)+v_p(a_1)+\ldots+v_p(a_n),\\
v_p(B_n)&=v_p(a_1)+v_p(a_2)+\ldots+v_p(a_n),
\end{split}
\end{equation}
or, equivalently,
\begin{equation}
\begin{split}
|A_n|_p&=|a_0|_p|a_1|_p\ldots|a_n|_p,\\
|B_n|_p&=|a_1|_p|a_2|_p\ldots|a_n|_p.
\end{split}
\end{equation}
The proof is done by induction on $n$. Using a similar argument, it is not hard to see that a sufficient condition for these equations to hold is having $v_p(a_{n+1}B_n)<v_p(B_{n-1})$ for all $n\in\mathbb{N}$.
\end{Remark}

The $p$--adic convergence strictly depends on the valuation of the partial quotients, regardless of the representative of $\mathbb{Z}/p\mathbb{Z}$ that are used. Therefore, the next proposition proves the convergence for both Ruban's and Browkin's continued fractions.

\begin{Proposition}[\cite{BI}]\label{Prop: ConvBr1}
Let an infinite sequence $a_0,a_1,\ldots\in \mathbb{Z}[\frac{1}{p}]$ such that for all $n\in\mathbb{N}$, $v_p(a_n)<0$. Then the continued fraction $[a_0,a_1,\ldots]$ is convergent to a $p$-adic number.
\end{Proposition}

\begin{Remark}
Proposition \ref{Prop: ConvBr1} and all the convergence results in this section contain conditions that hold for all $n$. However, the results do not change if these conditions hold for all $n\geq n_0\in\mathbb{N}$.
\end{Remark}

The proof of Proposition \ref{Prop: ConvBr1} is a simple induction which shows that, assuming  $v_p(a_n)<0$ for all $n\in\mathbb{N}$, then $\{v_p(B_n)\}_{n\in\mathbb{N}}$ is a strictly decreasing sequence. Therefore, $\{v_p(B_nB_{n+1})\}_{n\in\mathbb{N}}$ diverges to $-\infty$ and by the relation \eqref{Eq: conve} the continued fraction converges in $\mathbb{Q}_p$. However, requiring the sequence $\{v_p(B_nB_{n+1})\}_{n\in\mathbb{N}}$ strictly decreasing is equivalent to ask that $v_p(B_{n+1})<v_p(B_{n-1})$ for all $n\geq 1$. Therefore, the condition of Proposition \ref{Prop: ConvBr1} can be lightened by allowing partial quotient at even steps to have null valuation. In fact, the convergence of \textit{Browkin II} relies on the following proposition.

\begin{Proposition}[\cite{BII}]\label{Prop: ConvBr2}
Let an infinite sequence $a_0,a_1,\ldots\in \mathbb{Z}[\frac{1}{p}]$ such that for all $n\in\mathbb{N}$:
\begin{equation}
\begin{cases}
v_p(a_{2n})=0\\
v_p(a_{2n+1})<0.
\end{cases}
\end{equation}
Then the continued fraction $[a_0,a_1,\ldots]$ is convergent to a $p$-adic number.
\end{Proposition}
In both the proofs of Proposition \ref{Prop: ConvBr1} and Proposition \ref{Prop: ConvBr2}, the sequence $\{v_p(B_nB_{n+1})\}_{n\in\mathbb{N}}$ diverges to $-\infty$ because it is a strictly decreasing sequence of integers. In \cite{MRSII}, the strict decrease of this sequence has been effectively characterized in terms of valuations of the partial quotients. 

\begin{Theorem}[\cite{MRSII}]\label{Thm: ConvCond}
Let  $a_0,a_1,\ldots\in \mathbb{Z}\left[\frac{1}{p}\right]$ be an infinite sequence such that
\[v_p(a_{n}a_{n+1})<0,\]
for all $n > 0$.
Then the continued fraction $[a_0,a_1,\ldots]$ is convergent to a $p$--adic number.
\end{Theorem}

In the same paper, the authors explored some convergence conditions with partial quotients not satisfying the hypothesis of \ref{Thm: ConvCond}. They proved the following generalization of Proposition \ref{Prop: ConvBr2}.

\begin{Theorem}[\cite{MRSII}]\label{Thm: ConvBr3}
Let $a_0,a_1,\ldots \in \mathbb{Q}_p$ such that, for all $n\in\mathbb{N}$:
\[\begin{cases}
v_p(a_{3n+1})<0\\
v_p(a_{3n+2})=0\\
v_p(a_{3n+3})=0.
\end{cases}\]
If $v_p(a_{3n+3}a_{3n+2}+1)=0$ for all $n \in \mathbb{N}$, then, 
\[v_p(B_{3n-2})=v_p(B_{3n-1})=v_p(B_{3n})>v_p(B_{3n+1}).\]
\end{Theorem}

Theorem \ref{Thm: ConvBr3} allows to define Algorithm \eqref{Alg: Br3}, which satisfies both the convergence condition and the additional restriction on partial quotients. In \cite{MRSII}, Theorem \ref{Thm: ConvBr3} has been also generalized to an arbitrary number of steps by adding more hypotheses on the valuation of partial quotients.

\section{Continued fractions of rational numbers}\label{Sec: ratio}
In the field of real numbers, the rational numbers are characterized by finite continued fractions. This is a natural consequence of the finiteness of the Euclidean division algorithm. In fact, for all $a,b\in\mathbb{Z}$, $b\neq 0$, the simple continued fraction of $\frac{a}{b}$ is obtained by iterating
\begin{equation}
 \frac{a}{b}=\frac{bq+r}{b}=q+\frac{r}{b}=q+\frac{1}{\frac{b}{r}},  
\end{equation}
where $b=\lfloor a\rfloor$ and, hence, $|r|<|b|$. For $p$--adic versions of the Euclidean algorithm, see \cite{Err,Lager}.\\

In \cite{BUN}, Bundschuh proved that Schneider's continued fractions of rational numbers are not always finite, but they can also be periodic. In \cite{HW}, the authors gave a combinatorial characterization of some non-terminating expansions. Very recently, Pejkovic \cite{PEJ} proved the following effective criterion for determining whether Schneider's continued fraction of a rational number terminates.
\begin{Theorem}[\cite{PEJ}]
Let $\alpha=\frac{a}{b}\in\mathbb{Q}$. Schneider's continued fraction of $\alpha$ either terminates or a period is detected within $O(\log ^2 H(\alpha))$ steps, where $H(\alpha)=H(\frac{a}{b})=\max\{|a|,|b| \}$.
\end{Theorem}

For Ruban's continued fractions, Laohakosol \cite{LAO} proved the following, which is the analogue of Bundschuh's result.

\begin{Theorem}[\cite{LAO}]\label{Thm: Lao}
A $p$--adic number is rational if and only if either its Ruban's continued fraction terminates or it is eventually periodic with all partial quotients equal to $p-\frac{1}{p}$ from a certain point onward.
\end{Theorem}

\begin{Remark}
A $p$--adic number that has a finite Ruban's continued fraction is, by construction, rational and positive. Therefore, negative numbers can not have a finite Ruban's continued fractions. For example,

\begin{align*}
-p=(p-1)p+(p-1)p^2+\ldots=\cfrac{1}{\frac{p-1}{p}+(p-1)+\cfrac{1}{\frac{p-1}{p}+(p-1)+\ddots}},    
\end{align*}
that is, $-p=\left[\overline{p-\frac{1}{p}}\right]$.    
\end{Remark}

In the proof of Theorem \ref{Thm: Lao}, Laohakosol showed that all the continued fractions of rational numbers that are not finite, eventually have $-p$ as complete quotient. A proof of the characterization of rational numbers through this algorithm can be found also in \cite{WANG}. Ruban's continued fractions have been deepened in more details by Capuano, Veneziano and Zannier \cite{CVZ}. In particular, they provided an effective algorithm that determines in a finite number of steps whether the continued fraction of a rational number terminates or not. Moreover, in the non-terminating case, a negative complete quotients appears in the expansion of $\alpha=\frac{a}{b}$ in at most $\max\{2,\frac{\log b}{\log p}\}$ steps. Finally, they proved the following result that, giving a rational number $\alpha$, establishes how its Ruban's continued fraction changes when varying $p$.
\begin{Proposition}[\cite{CVZ}]
Let $\alpha\in\mathbb{Q}$. Then:
\begin{enumerate}
    \item[i)] If $\alpha<0$, then for every prime number $p$, Ruban's continued fraction of $\alpha$ does not terminate,
    \item[ii)] If $\alpha\geq 0$ and $\alpha\in\mathbb{Z}$, there are only finitely many prime numbers $p$ such that Ruban's continued fraction of $\alpha$ does not terminate,
    \item[iii)] If $\alpha\geq 0$ and $\alpha\not\in\mathbb{Z}$, there are only finitely many prime numbers $p$ such that Ruban's continued fraction of $\alpha$ terminates.
\end{enumerate}
\end{Proposition}

The first definition of $p$--adic continued fractions that terminate for every rational number is Browkin's first algorithm \eqref{Alg: Br1}. This is due to the choice of the representatives in $\{-\frac{p-1}{2},\ldots,\frac{p-1}{2}\}$ instead of $\{0,\ldots,p-1\}$. In this case, in fact, the Euclidean absolute value of the partial quotients satisfies the following inequality:
\begin{equation}\label{Eq: ineq}
|s(\alpha)|=\left|\sum\limits_{n=-r}^{0}c_np^n\right|\leq \frac{p-1}{2}\left|\sum\limits_{n=-r}^{0}p^n\right|<\frac{p}{2}.
\end{equation}
In order to give an idea of the proof of finiteness for Browkin's algorithm, let us write complete quotients of two consecutive steps
\[\alpha_n=\frac{A_n}{p^jB_n}, \ \ \alpha_{n+1}=\frac{A_{n+1}}{p^kB_{n+1}},\]
with $j,k\geq 1$, $(A_n,B_n)=(A_{n+1},B_{n+1})=1$ and $p\not|A_nB_nA_{n+1}B_{n+1}$. Then, by inequality \eqref{Eq: ineq} and the recursion
\[\alpha_{n+1}=\frac{1}{\alpha_n-a_n},\]
it is obtained, for all $n\in\mathbb{N}$,
\begin{equation}\label{Eq: finite}
|B_{n+1}|< \frac{|A_n|}{2}+\frac{|B_n|}{2}.
\end{equation}
The latter inequality allows to prove that $\{|A_n|+2|B_n|\}_{n\in\mathbb{N}}$ is a strictly decreasing sequence. Since it is a sequence of natural numbers, then it is finite. It means that there is only a finite number of complete quotients in Browkin's expansion of a rational number. Some results about the complexity of \textit{Browkin I} continued fractions for rational numbers can be found in \cite{BELE}. The finiteness of Browkin's second algorithm \eqref{Alg: Br2} is less straightforward, hence Browkin left it open in \cite{BII} as a conjecture. Most recently, Barbero, Cerruti and Murru \cite{BCMI} proved it by providing an inequality similar to \eqref{Eq: ineq} also for the second floor function $t$: for any $p$--adic number $\alpha=\sum\limits_{n=-r}^{+\infty}c_np^n$,
\begin{equation}\label{Eq: ineq2}
\left|t(\alpha)\right|=\left|\sum\limits_{n=-r}^{-1}c_np^n\right|<\frac{1}{2}.
\end{equation}
For the odd partial quotients $a_{2n+1}$, where the function $t$ is eventually adjusted with the sign function,
\begin{equation}\label{Eq: ineq3}
\left|a_{2n+1}\right|\leq 1-\frac{1}{p^l},
\end{equation}
where $l=-v_p(\alpha_{2n+1})$. Using \eqref{Eq: ineq} and \eqref{Eq: ineq3}, the authors provided two inequalities similar to \eqref{Eq: finite}, one for the odd steps and one for the even steps of the algorithm. This led to a strictly decreasing sequence of natural numbers similar to the one obtained by Browkin. Also the $p$--adic continued fractions algorithms defined in \cite{MR} and \cite{MRSII} terminate on rational inputs and the ideas for the proofs are similar. The proof of finiteness for the algorithm proposed in \cite{BCMII} is slightly different and it requires some further remarks on the absolute values of the denominators of the convergents. In fact, it has to handle with standard representatives in $\{0,\ldots,p-1\}$ and it is not possible to use inequalities \eqref{Eq: ineq} and \eqref{Eq: ineq2}.

\section{Periodicity properties}
\label{Sec: periodicity}
The continued fraction expansion of a real number is eventually periodic if and only it is a quadratic irrational number, by Lagrange's Theorem. Also for $p$--adic continued fractions, periodicity is often related to quadratic irrationals. In fact, a periodic continued fraction can be regarded as the root of a quadratic polynomial. However, we have seen in the previous section that periodic $p$--adic continued fractions can correspond also to rational numbers. On the other hand, none of the existent algorithms has been proved to produce a periodic continued fractions for all quadratic irrationals in $\mathbb{Q}_p$. Therefore, none of the two implications of Lagrange's Theorem is true in general for $p$--adic continued fractions. The latter is one of the most challenging open problems in this research area. This section is devoted to the main results regarding the periodicity of $p$--adic continued fractions up to the most recent developments.

\begin{Remark}
Periodic continued fractions, when not representing a rational number, converge to some irrational number $\alpha\in\mathbb{Q}_p$ that is quadratic over $\mathbb{Q}$, i.e. it is a root of an irreducible polynomial $f(x)\in\mathbb{Q}[x]$ of degree $2$. We denote by $\overline{\alpha}$ the conjugate of $\alpha$, that is the other root of $f(x)$ over $\mathbb{Q}$. Moreover, when we write $\sqrt{D}$ as a root of $x^2-D$ over $\mathbb{Q}_p$, it must be clear which of the two roots we are meaning. In the case of standard representatives in $\{0,\ldots,p-1\}$, we choose the root that in its $p$--adic expansion has the smaller first representative, while in the case of representatives in $\{-\frac{p-1}{2},\ldots,\frac{p-1}{2}\}$ we take the one with positive first representative. We denote the other root by $-\sqrt{D}$.
\end{Remark}

 Let $\alpha_0=\frac{P_0+\sqrt{D}}{Q_0}\in\mathbb{Q}_p$, where $D$ is a non-square integer that is a quadratic residue modulo $p$ and $P_0,Q_0\in\mathbb{Q}$. For all $n\in\mathbb{N}$, the complete quotients of the continued fraction expansion of $\alpha$ can be written as 
\begin{equation}\label{Eq: quadirr}
\alpha_n=\frac{P_n+\sqrt{D}}{Q_n},
\end{equation}
where $P_n,Q_n\in\mathbb{Q}$. Similarly as in the real framework, the sequences $P_n$ and $Q_n$, for all $n\in\mathbb{N}$ can be computed recursively, starting from $P_0$ and $Q_0$, by
\[P_{n+1}=a_nQ_n-P_n, \ \ \ \ \ \ \ Q_{n+1}=\frac{D-P^2_{n+1}}{Q_n},\]
where $a_n$ is the $n$-th partial quotient of $\alpha_0$.

\begin{Remark}\label{Rem: Lag}
A famous result of Galois \cite{GAL} states that the continued fraction expansion of a quadratic irrational $\alpha\in\mathbb{R}$ is purely periodic, i.e. periodic without a pre-periodic part, if and only if it is \textit{reduced}, that is $\alpha>1$ and its conjugate $-1<\overline{\alpha}<0$. One of technique to prove Lagrange's Theorem \cite{LAG} in $\mathbb{R}$ consists in showing that the expansion of any quadratic irrational eventually reaches a \textit{reduced} complete quotients, hence it starts to be periodic. Equivalently, it is possible to prove that there are only finitely many choices for the values of $P_n$ and $Q_n$ in \eqref{Eq: quadirr} for the expansions of a quadratic irrational. Therefore there is a finite number of complete quotients and the continued fractions must become periodic at some point.
\end{Remark}

One of the first works on the periodicity of $p$--adic continued fractions is due to Bundschuh \cite{BUN}. He studied Schneider's continued fractions and suggested that Lagrange's Theorem fails for this algorithm through some numerical computations, although he did not prove it. Some years later, de Weger \cite{DEWII} proved the following criterion for the non-periodicity of Schneider's continued fractions.
\begin{Proposition}[\cite{DEWII}]\label{Pro: Dewe}
Let $P_n$ and $Q_n$ as in \eqref{Eq: quadirr}. Then, if for some $n$ the signs of $P_n$ and $Q_n$ are different and $P_{n+1}^2>c$, Schneider's continued fraction of $\sqrt{c}\in\mathbb{Z}_p$ is not periodic. In particular, $\sqrt{c}\in\mathbb{Z}_p$ is never periodic for $c<0$.
\end{Proposition}
In \cite{DEWI}, the same author approached the periodicity of $p$--adic continued fractions from another point of view, i.e. by associating a sequence of \textit{approximation lattices} to every $p$--adic number. In the spirit of Lagrange's Theorem, de Weger proved that the sequence of approximation lattices attached to $\alpha\in\mathbb{Z}_p$ becomes periodic if and only if $\alpha$ is a quadratic irrational. However, this method is not effective for the construction of a periodic continued fraction for a given $p$--adic quadratic irrational. Few years later, Tilborghs \cite{TIL} determined an algorithm to detect the non-periodicity of Schneider's continued fraction in a finite number of steps. Becker \cite{BECK} then showed the following result on the length of the pre-periods.
\begin{Proposition}[\cite{BECK}]
Let $\alpha\in\mathbb{Z}_p$ be quadratic over $\mathbb{Q}$ with periodic Schneider's continued fraction. Then, if $p$ does not divide the discriminant of $\alpha$, the pre-period length is at most $1$.
\end{Proposition}
A fairly complete survey on the periodicity of Schneider's continued fractions is contained in \cite{VDP}. More recently, this algorithm has been studied also from other points of view and in more generality in \cite{NH, HJLN, HYZ}.\\

For Ruban's algorithm, using an argument similar to de Weger, Ooto \cite{OO} showed that not all quadratic irrationals in $\mathbb{Q}_p$ have a periodic Ruban's continued fraction. In particular, he proved the following result, giving some sufficient conditions to have a non-periodic $p$--adic continued fractions.

\begin{Proposition}[\cite{OO}]\label{Prop: Ooto}
Let $\alpha=\sqrt{D}\in\mathbb{Q}_p$ and the sequences $P_n,Q_n$ as in \eqref{Eq: quadirr} for some $n\in\mathbb{N}$, $P_nQ_n\leq 0$ and $P^2_{n+1}>D$, then the Ruban's continued fraction of $\alpha$ is not periodic.
\end{Proposition}

The idea of the proof of Proposition \ref{Prop: Ooto} is to show that, in these hypotheses, the sequence $P_n$ has strictly increasing absolute value from some point onward, hence it can not be periodic. From the latter proposition it easily follows that, for example, $\sqrt{D}$ can not have a periodic continued fraction if $D<0$. In fact, in this case, $P_0Q_0=0$ and $P_1\geq 0>D$ (in analogy to Proposition \ref{Pro: Dewe} for Schneider's algorithm).\\
In \cite{CVZ}, Capuano, Veneziano and Zannier made a more extensive analysis of  Ruban's algorithm. In particular, they provided an effective criterion for determining in a finite number of step whether the expansion of a quadratic irrational becomes periodic. They proved the following result.

\begin{Theorem}[\cite{CVZ}]\label{Thm: CVZ}
Let $\alpha\in\mathbb{Q}_p$ be a quadratic irrational over $\mathbb{Q}$. Then, Ruban's continued fraction of $\alpha$ is periodic if and only if there exists a unique real embedding $j:\mathbb{Q}(\alpha)\rightarrow \mathbb{R}$ such that the image of each complete quotient $\alpha_n$ under the map $j$ is positive. Moreover, there exists an effectively computable constant $N_{\alpha}$ with the property that, either exists $n\leq N_{\alpha}$ such that $\alpha_n$ does not have a positive real embedding, therefore the expansion is not periodic, or $\alpha_{n_1}=\alpha_{n_1}$ for $n_1<n_2\leq N_{\alpha}$, hence the expansion is periodic.
\end{Theorem}

The idea behind the proof of Theorem \ref{Thm: CVZ} is similar to the proof of Lagrange's Theorem pointed out in Remark \ref{Rem: Lag}. First of all, they proved an analogue of Galois' Theorem for classical continued fractions, finding a necessary condition for purely periodic continued fractions. This result is similar to Theorem \ref{Thm: PurePerBr1} for \textit{Browkin I}, since it is not affected by a different choice of representatives. Then, the authors used this necessary condition on the $p$--adic norm to prove that there are only finitely many quadratic irrational $\alpha_n$ of the form \eqref{Eq: quadirr} such that $\alpha_n$ has a purely periodic Ruban's continued fractions. Moreover, they showed that having two negative embeddings for some complete quotient $\alpha_n$, is a sufficient condition for the non-periodicity of the expansion. Finally, they  effectively computed a constant $N_{\alpha}$ such that, for $n\leq N_{\alpha}$, either one of the $\alpha_n$ has a purely periodic Ruban's continued fractions or it has two negative embeddings.\\

The situation is more complicated for Browkin's algorithms \eqref{Alg: Br1} and \eqref{Alg: Br2}, where the problem of deciding whether the $p$--adic continued fraction of a quadratic irrational is eventually periodic is still open. However, some numerical simulations suggest that this is not always the case (see, for example, \cite{BCMI,BCMII,BI, BII,CMT,MR}). As we have seen in Section \ref{Sec: ratio}, for these two algorithms finite continued fractions characterize rational numbers. Therefore, all periodic continued fractions correspond to irrational elements that are quadratic over $\mathbb{Q}$. The first systematic study of the periodicity of \textit{Browkin I} has been performed by Bedocchi in \cite{BEI, BEII, BEIII, BEIV}. The first results concern purely periodic continued fractions.

\begin{Theorem}[\cite{BEI}]\label{Thm: PurePerBr1}
Let $\alpha\in\mathbb{Q}_p$ having a periodic \textit{Browkin I} continued fraction expansion. Then the expansion is purely periodic if and only if
\[|\alpha|_p> 1, \quad |\overline{\alpha}|_p<1.\]
\end{Theorem}

This result is the analogue of Galois' Theorem for classical continued fractions and it is similar both in the statement and the proof. Moreover, he characterized the possible lengths of pre-periods for periodic \textit{Browkin I} continued fractions.

\begin{Proposition}[\cite{BEI}]
Let $D\in\mathbb{Z}$ not a perfect square, such that $\sqrt{D}\in\mathbb{Q}_p$. If the \textit{Browkin I} expansion of $\sqrt{D}$ is periodic, then the pre-period is
\begin{align*} 
\begin{cases}
2 \ \ \ \ \ & D \not\equiv 4 \mod 8 \ \text{when} \ p=2 \\
3 \ \ \ \ \ & otherwise.
\end{cases}
\end{align*}
\end{Proposition}

Then, in \cite{BEII, BEIII, BEIV}, Bedocchi focused on the possible lengths of the periods for square roots of integers. We collect all the results in the next proposition 

\begin{Proposition}[\cite{BEII, BEIII, BEIV}]\label{Prop: PropBedocchi}
For \textit{Browkin I}, the following statements are true.
\begin{enumerate}
    \item[i)] There are no periodic square roots of integers with period of length $1$;
    \item[ii)] For every odd integer $d$, there are only finitely many square roots of integers that are periodic with period of length $d$;
    \item[iii)] There exist infinitely many square roots of integers that are periodic with period of length $2$, $4$ and $6$.
\end{enumerate}   
\end{Proposition}

In view of the results of Proposition \ref{Prop: PropBedocchi}, Bedocchi conjectured that, for all $h$ even, there exist infinitely many square roots of integers that are periodic and have period length $h$. This conjecture is still open. However, for all the length that are powers of $2$, it has been proved by Capuano, Murru and Terracini in \cite{CMT}.

\begin{Theorem}[\cite{CMT}]\label{Thm: CMT}
For every $n,k\geq 1$, there are infinitely many $D\in\mathbb{Z}$, with $p\not | D$, such that the \textit{Browkin I} expansion of $p^k\sqrt{D}$ is periodic with period of length $2^n$. 
\end{Theorem}

The proof of Theorem \ref{Thm: CMT} exploits a new definition of a particular class of Browkin's continued fractions, that they call \textit{nice}, for which they proved the following theorem.

\begin{Theorem}[\cite{CMT}]\label{Thm: nice}
Let $[a_0,a_1,\ldots,a_{t-1}]$ be a \textit{nice} Browkin's continued fractions, for some $t\in\mathbb{N}$. Then there are infinitely many Browkin's partial quotients $a_t$ such that the continued fraction
\[[a_0,\overline{a_1,\ldots,a_{t-1},a_t,a_{t-1},\ldots,a_1,2a_0}]\]
converges to a quadratic irrational number of the form $\frac{1}{p^{e}\sqrt{D}}$, for some $D\in\mathbb{Z}$ not a perfect square.
\end{Theorem}

Basically, using Theorem \ref{Thm: nice}, starting with a \textit{nice} sequence $[a_0,a_1,\ldots,a_{t-1}]$ it is possible to provide an infinite family of integers $D$ such that the Browkin's continued fraction of $p^e\sqrt{D}$ is periodic with period $2t$. They left the following conjecture, together with other problems about \textit{nice} continued fractions.

\begin{Conjecture}
For every $t\geq 1$ there exist a \textit{nice} Browkin's continued fraction of length $t$, except when $t=1$ and $p=3$. 
\end{Conjecture}

By Theorem \ref{Thm: nice}, solving the latter conjecture implies the existence of infinitely many square roots of integers that have periodic \textit{Browkin I} expansion with period $2t$, for any $t$. Thus, it would give a positive answer to Bedocchi's conjecture. Moreover, in the same paper, the authors deepened other aspects of the periodicity of \textit{Browkin I}. The aim was to use some arguments similar to those of \cite{CVZ} in order to prove an effective criterion to predict the periodicity and the non-periodicity of \textit{Browkin I}. Unfortunately, Theorem \ref{Thm: CVZ} strongly depends on the fact that a periodic Ruban's continued fraction can always be embedded in the real numbers. This is not always true for Browkin's continued fractions, hence it is not a necessary condition for periodicity. On the other hand, having two negative embeddings is not a sufficient condition for non-periodicity (see \cite{CMT} for more details). This important difference between Ruban's and Browkin's continued fractions is again due to the choice of the (possibly negative) representatives in \textit{Browkin I}. For these reasons, there is not an ``easy" way to prove non-periodicity in \textit{Browkin I} up to now. In particular, no quadratic irrational has been proved to have non-periodic \textit{Browkin I} continued fraction, although it is largely believed to fail Lagrange's Theorem (see, for example, the experimental computations in \cite{MR}). The SageMath code that we have developed for \cite{MR} is publicly available\footnote{\href{https://github.com/giulianoromeont/p-adic-continued-fractions}{https://github.com/giulianoromeont/p-adic-continued-fractions}} and contains the implementation of the main algorithms presented in Section \ref{Sec: algo}. \\

In \cite{ABCM}, the authors followed another approach to find periodic $p$--adic continued fractions for all the square roots of integers in $\mathbb{Q}_p$. They did not focus on a specific algorithm and they used Rédei rational functions \cite{RED} to construct a periodic continued fraction converging to $\sqrt{D}$ simultaneously in the real and the $p$--adic field. It is possible to notice that, for any integer $z$, 
\begin{align*}
\sqrt{D}&=z+(\sqrt{D}-z)=z+\frac{1}{\frac{1}{\sqrt{D}-z}}=z+\frac{1}{\frac{\sqrt{D}+z}
{D-z^2}}=\\
&=z+\frac{1}{\frac{2z}{D-z^2}+\frac{\sqrt{D}-z}{D-z^2}}=z+\frac{1}{\frac{2z}{D-z^2}+\frac{1}{z+\sqrt{D}}},
\end{align*}
hence
\begin{equation}\label{Eq: Red}
\sqrt{D}=\left[z,\overline{\frac{2z}{d-z^2},2z} \right],
\end{equation}
even if it does not coincide with its expansion with the standard algorithm (see \cite{PER2}). Rédei rational functions are known to converge to $\sqrt{D}$ in $\mathbb{R}$ (see, for example, \cite{DICK}). In \cite{ABCM}, the authors proved that Rédei functions converge to $\sqrt{D}$ also in $\mathbb{Q}_p$, so that \eqref{Eq: Red} is true also for $p$--adic numbers. This construction has been generalized in \cite{BCMI} in order to manage any quadratic irrational, not only square roots of integers. The authors proposed a generalization of Rédei rational functions to provide a periodic continued fraction for $\alpha$ root of the polynomial $x^2+hx-d$, where $h,d\in\mathbb{Z}$. The expansion is
\begin{equation}\label{Eq: exp}
\alpha=\left[z,-\overline{\frac{h+2z}{z^2+hz-d},h+2z}  \right],
\end{equation}
that yields a periodic representation for $\alpha$ both in $\mathbb{R}$ and $\mathbb{Q}_p$. This result solves the problem of expressing every $p$--adic quadratic irrational number as a periodic continued fraction, i.e. through a periodic sequence of rational numbers. However, it is not equivalent to Lagrange's Theorem, since the expansion is not obtained by a specific algorithm. In fact, in order to find the expansion \eqref{Eq: exp} for $\alpha$, we need from the begin the knowledge that it is a quadratic irrational number and, in addition, of its characteristic polynomial over $\mathbb{Q}$. In the same paper, the authors characterized the cases when \textit{Browkin II} provides a periodic expansion of the form \eqref{Eq: exp}.

\begin{Theorem}[\cite{BCMI}]
Given $\alpha\not\in\mathbb{Q}$ root of the polynomial $x^2+hx-d$, with $h,d\in\mathbb{Z}$, \textit{Browkin II} produces the $p$--adic continued fraction
\begin{equation}
\alpha=\left[z,-\overline{\frac{h+2z}{p},h+2z}  \right],
\end{equation}
if and only if
\[1\leq|z|\leq\frac{p-1}{2}, \ \ \ 1\leq|h+2z|\leq\frac{p-1}{2},\]
for $z$ such that $z^2+hz-d=p$.
\end{Theorem}

The authors left open the problem of finding an actual algorithm that, for all quadratic irrationals in $\mathbb{Q}_p$, provides periodic representations of the form \eqref{Eq: exp}. Recently, the properties of periodicity of Browkin's second algorithm \eqref{Alg: Br2} have been considered. In \cite{MRSI}, an approach similar to Bedocchi in \cite{BEI,BEII,BEIII,BEIV} is employed to characterize purely periodic continued fractions and to study the lengths of pre-periods and periods for periodic expansion. The situation is more complicated than \textit{Browkin I}, because of the presence of the sign function in \eqref{Alg: Br2} and it is possible to obtain only some partial results. In fact, following an argument similar to the proof of Theorem \ref{Thm: PurePerBr1}, the condition on the $p$--adic absolute values of $\alpha$ and its conjugate are only necessary.

\begin{Theorem}[\cite{BEI}]\label{Thm: Bedo1}
If $\alpha\in\mathbb{Q}_p$ has a purely periodic continued fraction expansion
\[\alpha=[\overline{a_0,\ldots,a_{k-1}}],\]
with \textit{Browkin II}, then
\[ |\alpha|_p=1, \ \ |\overline{\alpha}|_p<1.\]
Conversely, if $\alpha$ has a periodic \textit{Browkin II} expansion
\[\alpha=[a_0,a_1,\ldots,a_{h-1},\overline{a_h,\ldots,a_{h+k-1}}],\]
and $|\alpha|_p=1,|\overline{\alpha}|_p<1$, then the preperiod length can not be odd.
\end{Theorem}
Theorem \ref{Thm: Bedo1} can not be improved, since there are counterexamples of periodic continued fractions for elements that satisfy hypotheses, but are not purely periodic. Moreover, in the spirit of Bedocchi's Proposition \ref{Prop: PropBedocchi}, some results have obtained on the possible pre-periods and periods of \textit{Browkin II} expansions.

\begin{Proposition}[\cite{MRSI}] The following results hold for \textit{Browkin II}:
\begin{enumerate}
\item[i)] If $\sqrt{D}\in\mathbb{Q}_p$ has a periodic continued fraction, then the pre-period has length either $1$ or even;
\item[ii)] There are infinitely many $\sqrt{D} \in \mathbb{Q}_p$ having a periodic $p$--adic expansion with period $4$.
\end{enumerate}
\end{Proposition}

Supported by experimental computations, the authors left the following conjecture, which is the exact analogue of Bedocchi's conjecture for Browkin's second algorithm.

\begin{Conjecture}[\cite{MRSI}]
For all even $h\in\mathbb{Z}$, there exist infinitely many $\sqrt{D}$, $D\in\mathbb{Z}$ not perfect square, such that \textit{Browkin II} continued fraction of $\sqrt{D}$ is periodic with period of length $h$.
\end{Conjecture}

In \cite{MR}, the authors proved some further results on the periodicity \textit{Browkin II}, underlining some issues arising from the unpredictability of the appearance of the sign function in the algorithm. Algorithm \eqref{Alg: new}, introduced in \cite{MR}, improves \textit{Browkin II} from both a theoretical and an experimental point of view.

\begin{Theorem}[\cite{MR}]\label{Thm: newalg}
The following results hold for Algorithm \eqref{Alg: new}:
\begin{enumerate}
\item[i)] If $\alpha\in\mathbb{Q}_p$ has a periodic continued fraction, then it is purely periodic if and only if $|\alpha|_p\geq 1$ and $|\overline{\alpha}|_p<1$;
\item[ii)] If $\sqrt{D}$ has a periodic continued fraction, then the pre-period is $1$ for $v_p(\sqrt{D})\leq 0$ and $2$ for $v_p(\sqrt{D})> 0$
\end{enumerate}
\end{Theorem}

Theorem \ref{Thm: newalg} establishes that an analogue of Galois' Theorem is true for Algorithm \eqref{Alg: new}, unlike Theorem \ref{Thm: Bedo1} for \textit{Browkin II}. Moreover, the second part is very similar to what happens in $\mathbb{R}$, where the pre-periods of square roots of integers is always $1$.

\section{Quality of the approximation}\label{Sec: approxi}
Classical continued fractions are one of the most powerful tools to produce rational approximations of real numbers. In fact, the convergents provide at each step the best approximation of real numbers, in the sense of the following proposition.
\begin{Proposition}\label{Pro: bestapprox}
Let $\alpha\in\mathbb{R}$ and $\{\frac{A_n}{B_n}\}_{n\in\mathbb{N}}$ the sequence of convergents of $\alpha$. Then
\[|B_n\alpha-A_n|\leq |B\alpha-A|,\]
for all $n\in\mathbb{N}$ and $A,B$  integers such that $0<B\leq B_n$.
\end{Proposition}

\begin{Remark}
In other words, Proposition \ref{Pro: bestapprox} tells us that there is no rational number with denominator minor than $B_n$ that approximates $\alpha$ better than $\frac{A_n}{B_n}$ and, fixing $B_n$ as denominator, the numerator $A_n$ is the best possible. Therefore, the approximations provided by the convergents of a continued fraction are the \textit{best approximations} without increasing the denominator (of course, with a bigger denominator, it is possible to obtain better approximations). 
\end{Remark}

In this section we analyze the quality of the approximations provided by continued fractions in $\mathbb{Q}_p$. First of all, let us give a more precise shape to the approximations given by the convergents of a $p$--adic continued fraction. Let $\alpha\in\mathbb{Q}_p$ such that $\alpha = [a_0, a_1, \ldots]$. We can write $\alpha=[a_0,a_1,\ldots,a_n,\alpha_{n+1}]$, so that, by \eqref{Eq: Recursions},
\[\alpha=\frac{\alpha_{n+1}A_n+A_{n-1}}{\alpha_{n+1}B_n+B_{n-1}}.\]
Therefore, the difference between $\alpha$ and its $n$-th convergent is
\[\alpha-\frac{A_n}{B_n}=\frac{\alpha_{n+1}A_n+A_{n-1}}{\alpha_{n+1}B_n+B_{n-1}}-\frac{A_n}{B_n}=\frac{(-1)^n}{(\alpha_{n+1}B_n+B_{n-1})B_n}.\]
Since $v_p(\alpha_{n+1})=v_p(a_{n+1})$, then $v_p(\alpha_{n+1}B_n)<v_p(B_{n-1})$, hence
\[v_p(\alpha_{n+1}B_n+B_{n-1})=v_p(\alpha_{n+1}B_n)=v_p(a_{n+1}B_n)=v_p(B_{n+1}).\] The valuation of the difference is then
\begin{equation}\label{Eq: val2}
v_p\left(\alpha-\frac{A_n}{B_n}\right)=-v_p(B_{n}B_{n+1}),
\end{equation}
that is, 
\begin{equation}\label{Eq: val}
\left|\alpha-\frac{A_n}{B_n}\right|_p=p^{v_p(B_nB_{n+1})}.
\end{equation}
By Remark \ref{Rem: ConveBr}, we know that, for all $n\in\mathbb{N}$,
\[v_p(B_n)=v_p(a_0)+v_p(a_1)+\ldots+v_p(a_n).\]
Therefore, by \eqref{Eq: val2} and \eqref{Eq: val}, the quality of approximation of $\alpha$ depends on the size of the valuation of the denominators $v_p(B_n)$. In fact, the algorithm that better approximates a $p$--adic number $\alpha$ is the one where $v_p(B_n)$ becomes ``very negative" rapidly.
\begin{Remark}\label{Rem: appbrow}
By Proposition \ref{Prop: ConvBr1} and Proposition \ref{Prop: ConvBr2}, we know that \textit{Browkin I} decreases the valuation of $B_n$ at each step, while \textit{Browkin II} only at odd steps. Therefore, at each step $n\in\mathbb{N}$, we expect \textit{Browkin I} to produce better approximations than \textit{Browkin II}.
\end{Remark}
Some computational results on the quality of approximations for square roots of integers by means of Algorithms \eqref{Alg: Br1}, \eqref{Alg: Br2} and \eqref{Alg: new} can be found in \cite{MR} and they empirically confirm the insight of Remark \ref{Rem: appbrow}.\\

In \cite{MAH}, Mahler attempted to obtain a sequence of best approximations for a $p$--adic number $\alpha$, similarly to the sequence of convergents for real continued fractions. Let $\Phi(X,Y)$ be a reduced positive definite quadratic form of determinant $-1$ (this is in analogy with the determinant of the matrices appearing in \eqref{Eq: matrix}). Mahler defined an algorithm to provide a sequence of pairs of integers $\{(A_n,B_n)\}_{n\in\mathbb{N}}$, such that, for all $n\in\mathbb{N}$,
\begin{equation}\label{Eq: bestap}
    |A_n-B_n\alpha|_p\leq\frac{1}{p^n},
\end{equation}
and
\[|A_n-B_n\alpha|_p<|A-B\alpha|_p,\]
for all the other pairs $(A,B)$ with $0<\Phi(A,B)<\Phi(A_n,B_n)$. Notice that this result is very similar to Proposition \ref{Pro: bestapprox} for real continued fractions, i.e. it tries to reproduce the idea of best approximations also in the $p$--adic framework. In \cite{DEA}, Deanin focused on Mahler's approach and studied the periodicity of the sequence $(A_n,B_n)$ by varying the function $\Phi$. As already mentioned in Section \ref{Sec: periodicity}, de Weger \cite{DEWI} studied rational approximations in $\mathbb{Q}_p$ by associating to each $p$--adic number a sequence of \textit{approximation lattices}. For all $n\in\mathbb{N}$, the $n$-th approximation lattice of $\alpha\in\mathbb{Q}_p$ is
\[\Gamma_n=\{(A,B)\in\mathbb{Z}^2 : |B\alpha-A|_p\leq p^{-n}\},\]
i.e. the lattice containing all the pairs of integers satisfying \eqref{Eq: bestap}. In this context, de Weger was able to prove a form of Hurwitz's Theorem \cite{HUR}, together with other known results of Diophantine approximation, over $\mathbb{Q}_p$. In \cite{BUG,BUG2,BBDO,BAD}, Schneider's $p$--adic continued fractions are employed for the study of simultaneous uniform approximations for a $p$--adic number and its powers. These problems are similar to those addressed in the real framework in \cite{BUG0,LAU,BAD} (see also \cite{BUGEAUD} for a survey on the topic). Simultaneous approximations of two $p$--adic numbers have been studied also in \cite{MT3} using multidimensional continued fractions. In fact, by definition, continued fractions in two dimensions provide simultaneous approximations by means of two sequences of integers (more details are given in Section \ref{Sec: multidim}, which is devoted to the study of multidimensional $p$-adic continued fractions).

\section{Transcendental $p$--adic continued fractions}\label{Sec: transce}
Continued fractions have been proved to be very efficient for the construction of transcendental numbers. Over the real numbers, Khinchin \cite{KHI} addressed the problem of continued fraction expansions of algebraic irrational numbers which are not quadratic over $\mathbb{Q}$. Liouville \cite{LIO} was the first to construct transcendental numbers by using continued fractions with unbounded partial quotients. Transcendental continued fractions with bounded partial quotients have been defined by Maillet \cite{MAI} and Baker \cite{BAK}. An extensive analysis has been lately performed in \cite{AB, AB2, AB3, ABD}, where the authors showed several transcendence criteria and provided many examples of families of transcendental continued fractions. One of the most powerful tool is the use of Roth's Theorem \cite{ROTH} and its subsequent improvements, in particular the  Schmidt's Subspace Theorem \cite{SCHM}. Let us recall Roth's Theorem in one of its simplest forms.
\begin{Theorem}[Roth's Theorem]
Let $\alpha$ be an algebraic irrational over $\mathbb{Q}$ and let $\epsilon >0$. Then there exist only finitely many $(a,b)\in\mathbb{Z}^2$, $b>0$, such that
\[\left|\alpha-\frac{a}{b}\right|<\frac{1}{b^{2+\epsilon}}.\]  
\end{Theorem}
Roth's Theorem tells us that an algebraic irrational number admits only finitely many ``good" approximations. Therefore, it is possible to prove the transcendence of a real number by constructing infinitely many ``good enough" rational approximations. In \cite{SCHL}, Schlickewei proved the following $p$--adic version of the Subspace Theorem.
\begin{Theorem}[\cite{SCHL}]\label{Thm: Schlik}
Let $\mathbf{x}=(x_1,\ldots,x_n)\in\mathbb{Q}^n$ and let us define
\[\Vert \mathbf{x} \Vert_{\infty}=\max \{|x_i|,i=1,\ldots,n\}, \ \ \Vert \mathbf{x} \Vert_{p}=\max \{|x_i|_p,i=1,\ldots,n\}.\]
Let $L_{1,\infty}(\mathbf{x}),\ldots,L_{n,\infty}(\mathbf{x})$ and $L_{1,p}(\mathbf{x}),\ldots,L_{n,p}(\mathbf{x})$ be, respectively, real and $p$--adic linearly independent linear forms with algebraic coefficients and let $\epsilon>0$. Then, the solutions $\mathbf{x}\in\mathbb{Z}^n$ of
\[\prod\limits_{i=1}^n(|L_{i,\infty}(\mathbf{x})|\cdot |L_{1,p}(\mathbf{x})|_p)\leq \frac{1}{\Vert \mathbf{x} \Vert_{\infty}^{\epsilon}},\]
lie in finitely many proper subspaces of $\mathbb{Q}^n$.
\end{Theorem}

Theorem \ref{Thm: Schlik} has been employed in \cite{BEH} to prove a transcendence result for the $p$--adic Thue-Morse continued fractions, inspired by Queffélec's result over the real numbers \cite{QUE}. In \cite{OO}, Ooto studied transcendental Ruban's $p$--adic continued fractions, performing an analysis similar to Baker \cite{BAK}. Baker's result for continued fractions in $\mathbb{R}$ relies on the definition of quasi-periodic continued fractions, i.e. continued fractions where blocks of partial quotients repeat several times.
\begin{Definition}\label{Def: quasiper}
Let $\{n_i\}_{i\in\mathbb{N}}$, $\{\lambda_i\}_{i\in\mathbb{N}}$ and $\{k_i\}_{i\in\mathbb{N}}$ be sequences of positive integers. If $n_{i+1}\geq n_i+\lambda_ik_i$ and $a_{m+k_i}=a_m$ for all $m=n_i,\ldots,n_i+(\lambda_i-1)k_i-1$, then the continued fraction $[a_0,a_1,\ldots]$ is called \textit{quasi-periodic}. 
\end{Definition}

The main result of \cite{OO} is the following transcendence criterion for quasi-periodic Ruban's continued fractions.

\begin{Theorem}[\cite{OO}]
Let $[a_0,a_1,\ldots]$ be a quasi-periodic Ruban's continued fraction which is not periodic and let $A\geq p$. Let  $a_{n_i}=\ldots=a_{n_i+k_i-1}=p-p^{-1}$ for infinitely many $i$, where the $n_i$ and $k_i$ are as in Definition \ref{Def: quasiper}. Then, if
\[\liminf\limits_{i\rightarrow +\infty} \frac{\lambda_i}{n_i}>B,\]
where $B=\frac{2\log A}{\log p}-1$, then Ruban's continued fraction $[a_0,a_1,\ldots]$ is transcendental.
\end{Theorem}

Other results on the transcendence of Ruban's continued fractions are obtained in \cite{AMM1,AMM2}, exploiting Theorem \ref{Thm: Schlik} and using a method similar to \cite{ABD} for the real case. A study of the transcendence Ruban's continued fraction has been performed also in \cite{DL, LV, WANG, DEWII}, using unbounded partial quotients. Finally, very recently, a transcendence criterion has been proved in \cite{LMS} for \textit{Browkin I}. We collected the main results in the following theorem, which contains the $p$--adic analogues of the transcendence criteria proved by Baker in \cite{BAK} and by Adamczewski and Bugeaud in \cite{AB3}.

\begin{Theorem}[\cite{LMS}]
Let $\alpha=[0,a_1,a_2,\ldots]$ be a non-periodic \textit{Browkin I} $p$--adic continued fraction. Let $\{n_i\}_{i\in\mathbb{N}}$, $\{\lambda_i\}_{i\in\mathbb{N}}$ and $\{k_i\}_{i\in\mathbb{N}}$ be such that  $n_{i+1}\geq n_i+\lambda_ik_i$, with $k_i$ bounded, and suppose there exist $\delta>0$ and $B$ such that $\lambda_i>(B+\delta)n_i$ for all $i$. Assume that one of the following holds:
\begin{enumerate}
    \item[i)] the sequence $\{|a_n|_p\}_{n\geq 1}$ is bounded and there are infinitely many subsequences
    \[a_{n_i}=a_{n_i+1}=\ldots =a_{n_i+k_i-1}=p^{-1},\]
    \item[ii)] the sequence $\{|a_n|_p\}_{n\geq 1}$ begins with arbitrarily long palindromes,
\end{enumerate}
where the $n_i$ and $k_i$ are as in Definition \ref{Def: quasiper}. Then $\alpha$ is transcendental.
\end{Theorem}

\section{Multidimensional $p$--adic continued fractions}\label{Sec: multidim}
The theory of continued fractions has been extended in several senses. One of the most famous generalization is given by multidimensional continued fractions. Motivated by a question of Hermite \cite{HER}, Jacobi \cite{JAC} introduced a continued fraction algorithm with the purpose of providing periodic representations for cubic irrationals, in the spirit of Lagrange's Theorem for quadratic irrationals. Perron \cite{PER} generalized Jacobi's approach in order to handle algebraic irrationalities of any degree. Hermite's problem for multidimensional continued fractions is still open. We refer the interested readers to \cite{BER, MUL} for the general theory of multidimensional continued fractions and to \cite{OK1, OK2} for some of the most recent developments on the Hermite's problem. One of the first work on multidimensional continued fractions in the field of $p$--adic numbers is due to Ruban \cite{RUB2, RUB3,RUB4}, who studied some of the ergodic properties of Jacobi-Perron algorithm in $\mathbb{Q}_p$. In the very last years, other multidimensional continued fractions algorithms have been defined in $\mathbb{Q}_p$. In particular, in \cite{MT1} and \cite{STY} the authors introduced two different multidimensional $p$--adic continued fraction algorithms trying to generalize, respectively, Browkin's and Schneider's approach (see also \cite{RY}). In \cite{MT2}, the authors deepened the study of the algorithm from \cite{MT1}, by exploring its properties of finiteness.\\

The classical Jacobi-Perron algorithm for multidimensional continued fraction works on a $m$-tuple of real numbers $(\alpha_0^{(1)},\ldots,\alpha_0^{(m)})$ and provides their representations through an $m$-tuple of integer sequences $(\{a_n^{(1)}\}_{n\in\mathbb{N}},\ldots,\{a_n^{(m)}\}_{n\in\mathbb{N}})$, generated by the following algorithm, for all $n\in\mathbb{N}$:
\begin{align*}
\begin{cases}
a_n^{(i)}=\lfloor \alpha_n^{(i)} \rfloor, \ \ i=1,\ldots,m\\
\alpha_{n+1}^{(i)}=\frac{1}{\alpha_n^{(m)}-a_n^{(m)}}\\
\alpha_{n+1}^{(i)}=\frac{\alpha_n^{(i-1)}-a_n^{(i-1)}}{\alpha_n^{(m)}-a_n^{(m)}}, \ \ i=2,\ldots,m.
\end{cases}
\end{align*}
Therefore, the Jacobi-Perron algorithm express every real number of the $m$-tuple $(\alpha_0^{(1)},\ldots,\alpha_0^{(m)})$ as a continued fraction by
\begin{align*} 
\begin{cases}
\alpha_{n}^{(i-1)}=a_{n}^{(i-1)}+\frac{\alpha_{n+1}^{(i)}}{\alpha_{n+1}^{(1)}}\\
\alpha_{n}^{(m)}=a_{n}^{(m)}+\frac{1}{\alpha_{n+1}^{(1)}}.
\end{cases}
\end{align*}
Also for multidimensional continued fractions it is possible to introduce the  convergents $\frac{A_n^{(i)}}{A_n^{(m+1)}}$, where the quantities $A_n^{(i)}$ are defined as follows, for all $n\in\mathbb{N}$ and for all $i=1,\ldots,m$:
\begin{align*}
\begin{cases}
A_{-j}^{(i)}=\delta_{ij}\\
A_0=a_0^{(i)}\\
A_n^{(i)}=\sum\limits_{j=1}^{m+1}a_n^{(j)}A_{n-j}^{i},
\end{cases}
\end{align*}
with the convention that $a_n^{(m+1)}=1$ for all $n\in\mathbb{N}$. The $p$--adic Jacobi-Perron algorithm presented in \cite{MT1} uses Browkin's $s$ function \eqref{Eq: sfunc} in place of the floor function of the classical algorithm. Hence, it works as follows, on an input $m$-tuple of $p$--adic numbers $(\alpha_0^{(1)},\ldots,\alpha_0^{(m)})$:

\begin{align}\label{Alg: padicJP}
\begin{cases}
a_n^{(i)}=s(\alpha_n^{(i)}), \ \ i=1,\ldots,m\\
\alpha_{n+1}^{(i)}=\frac{1}{\alpha_n^{(m)}-a_n^{(m)}},\\
\alpha_{n+1}^{(i)}=\frac{\alpha_n^{(i-1)}-a_n^{(i-1)}}{\alpha_n^{(m)}-a_n^{(m)}}, \ \ i=2,\ldots,m.
\end{cases}
\end{align}

The first results in \cite{MT1} concern the $p$--adic convergence of a generic multidimensional continued fractions. Some conditions on the $p$--adic absolute value similar to those seen in Section \ref{Sec: convergence} are necessary in order to guarantee convergence.

\begin{Proposition}[\cite{MT1}]
Let $(\{a_n^{(1)}\}_{n\in\mathbb{N}},\ldots,\{a_n^{(m)}\}_{n\in\mathbb{N}})$ be the sequences of partial quotients of a multidimensional $p$--adic continued fraction with, for all $n\in\mathbb{N}$,
\begin{align*}
\begin{cases}
v_p(a_n^{(1)})\leq 0\\
v_p(a_n^{(1)})<v_p(a_n^{(i)}), \ \ \ i=2,\ldots,m+1.
\end{cases}
\end{align*}
Then, for all $i=1,\ldots,m$, each of the sequences $\left(\frac{A_n^{(i)}}{A_n^{(m+1)}}\right)_{n\in\mathbb{N}}$ converges to a $p$--adic number.
\end{Proposition}
When the $p$--adic Jacobi-Perron converges to an $m$-tuple of $p$--adic numbers
\[(\alpha_0^{(1)},\ldots,\alpha_0^{(m)}),\]
then it is also strong convergent, that is,
\[\lim_{n\rightarrow +\infty}|A_n^{(i)}-\alpha_0^{(i)}A_n^{(m+1)}|_p=0,\]
for all $i=1,\ldots,m$. This is not the case in $\mathbb{R}$ where strong convergence is not guaranteed in general. Unlike classical continued fractions in one dimension, finite multidimensional continued fractions do not characterize rational numbers.

\begin{Proposition}[\cite{MT1}]
If the $p$--adic Jacobi-Perron algorithm \eqref{Alg: padicJP} stops in a finite number of steps on input $(\alpha_0^{(1)},\ldots,\alpha_0^{(m)})\in\mathbb{Q}_p^m$, then $1,\alpha_0^{(1)},\ldots,\alpha_0^{(m)}$ are $\mathbb{Q}$-linearly dependent.
\end{Proposition}

\begin{Theorem}[\cite{MT1}]\label{Thm: finitejp}
If $(\alpha_0^{(1)},\ldots,\alpha_0^{(m)})\in\mathbb{Q}^m$, then the $p$--adic Jacobi-Perron algorithm \eqref{Alg: padicJP} stops in a finite number of steps.
\end{Theorem}
The conditions of both the previous results are only sufficient and not necessary. In fact, there are $\mathbb{Q}$-linearly dependent inputs for which the algorithm is not finite but periodic, and also $\mathbb{Q}$-linearly dependent inputs having finite continued fraction although they are not all rationals (see \cite{MT1} for examples). Moreover, in \cite{MT2} the same authors characterized some of the $\mathbb{Q}$-linearly dependent inputs giving rise to finite continued fractions.

\begin{Theorem}[\cite{MT2}]
Let $1,\alpha_0^{(1)},\ldots,\alpha_0^{(m)}$ be $\mathbb{Q}$-linearly dependent $p$--adic numbers, with
\[v_p(a_n^{(j)})-v_p(a_n^{(1)})\geq j-1,\]
for $j=3,\ldots,m+1$ and $n$ sufficiently large. Then the Jacobi-Perron algorithm stops in finitely many steps on input $(\alpha_0^{(1)},\ldots,\alpha_0^{(m)})$.
\end{Theorem}

Therefore, the study of the expansion of $\mathbb{Q}$-linearly dependent $p$--adic numbers by means of the $p$--adic Jacobi-Perron algorithm is not complete. However, it has not been observed an $m$-tuple of $\mathbb{Q}$-linearly dependent $p$--adic numbers for which the continued fraction is infinite and not periodic. Therefore, in \cite{MT1} the authors left the following (still open) conjecture.

\begin{Conjecture}
Let $1,\alpha_0^{(1)},\ldots,\alpha_0^{(m)}$ be $\mathbb{Q}$-linearly dependent $p$--adic numbers. Then the Jacobi-Perron algorithm on input $(\alpha_0^{(1)},\ldots,\alpha_0^{(m)})$ is either finite or periodic.
\end{Conjecture}

For the periodicity of the $p$--adic Jacobi-Perron algorithm, the following result holds, similarly to multidimensional continued fractions in $\mathbb{R}$.

\begin{Theorem}[\cite{MT2}]
A periodic $m$-dimensional continued fractions obtained by the $p$--adic Jacobi-Perron algorithm \eqref{Alg: padicJP}, converges to an $m$-tuple of algebraic irrationalities of degree less or equal than $m+1$.
\end{Theorem}

In \cite{MT3}, the same authors explored the properties of approximation of the Jacobi-Perron algorithm in two dimensions, i.e. for the simultaneous approximation of two $p$--adic numbers. They studied the rate of convergence and gave some precise results for length of the finite continued fraction of two rational numbers, following the technique of the proof of finiteness of Theorem \ref{Thm: finitejp}. Moreover, they focused on the study of algebraic dependence of two $p$--adic numbers. In particular, they proved the following sufficient condition to ensure the finiteness of the $p$--adic Jacobi-Perron algorithm on $\mathbb{Q}$-linearly dependent inputs.
\begin{Theorem}[\cite{MT3}]
Let $\alpha,\beta\in\mathbb{Q}_p$ and let $\left(\frac{A_n}{C_n},\frac{B_n}{C_n}\right)$ be, for all $n\in\mathbb{N}$, the sequences of convergents of the multidimensional continued fraction representing $(\alpha,\beta)$. Moreover, let us define
\[M_n=\max\{|A_n|_{\infty},|B_n|_{\infty},|C_n|_{\infty}\}, \ \ \ U_n=\max\left\{\left|\alpha-\frac{A_n}{C_n}\right|,\left|\alpha-\frac{B_n}{C_n}\right|\right\}.\]
Then, if
\[\lim\limits_{n\rightarrow +\infty}U_n\cdot M_n=0,\]
either $1,\alpha,\beta$ are $\mathbb{Q}$-linearly independent or the expansion of $(\alpha,\beta)$ is finite.
\end{Theorem}
Finally, they found some conditions to define a $p$--adic multidimensional continued fraction converging to algebraically independent numbers.

\section*{Acknowledgments}
The author is a member of GNSAGA of INdAM.

\end{document}